\documentclass[12pt]{article}
\usepackage{amsmath,amssymb}
\usepackage[mathscr]{eucal}
\newtheorem{Thm}{Theorem}[section]
\newtheorem{Def}[Thm]{Definition}

\newtheorem{Prop}[Thm]{Proposition}

\newcommand{\Section}[1]{\setcounter{equation}{0}\section{#1}}

\newcommand{\RR}{{\mathbb R}}

\newcommand{\NN}{{\mathbb N}}
\newcommand{\refeq}[1]{{\rm(\ref{#1})}}
\newcommand{\proof}{\noindent{\bf Proof:}\quad}
\newcommand{\qed}{\hfill{Q.E.D.}\vspace{2mm}}
\newcommand{\vs}{\vspace{2mm}}
\newcommand{\Vs}[1]{\rule{0mm}{#1mm}}
\newcommand{\vsn}{\vspace{-1.6mm}}
\newcommand{\vsnp}{\vspace{-4mm}}
\newcommand{\vsnn}{\vspace{-1mm}}
\newcommand{\vsnpp}{\vspace{-4mm}}
\newcommand{\tens}{\otimes}
\newcommand{\id}{{\rm id}}
\newcommand{\Mor}{{\rm Mor}}
\newcommand{\is}[2]{\left(#1\,\vline\,#2\right)}
\newcommand{\its}[3]{\left(#1\,\vline\,#2\,\vline\,#3\right)}
\newcommand{\nc}{\raisebox{.5ex}{\mbox{\rm\tiny norm closure}}}
\newcommand{\Del}{\Delta}
\newcommand{\ka}{\kappa}
\newcommand{\Wtil}{\widetilde{W}}
\newcommand{\Qhat}{\widehat{Q}}
\newcommand{\Hbar}{\overline{H}}

\newcommand{\xbar}{\overline{x}}
\newcommand{\zbar}{\overline{z}}
\newcommand{\Ahat}{\widehat{A}}
\newcommand{\What}{\widehat{W}}
\newcommand{\Whattil}{\widetilde{\widehat{W}}}
\newcommand{\Tau}{{\mathscr T}}
\newcommand{\Rm}{R_{\mbox{\rm\tiny M}}}
\newcommand{\Wm}{W_{\mbox{\rm\tiny M}}}
\newcommand{\Wmtil}{\widetilde{W}_{\mbox{\rm\tiny M}}}
\newcommand{\Qm}{Q_{\mbox{\rm\tiny M}}}
\newcommand{\Am}{A_{\mbox{\rm\tiny M}}}
\newcommand{\am}{a_{\mbox{\rm\tiny M}}}
\newcommand{\Amhat}{\widehat{A}_{\mbox{\rm\tiny M}}}
\newcommand{\kam}{\kappa_{\mbox{\rm\tiny M}}}
\newcommand{\taum}{{\tau_{\mbox{\rm\tiny M}}}}
\newcommand{\Hm}{H_{\mbox{\rm\tiny M}}}
\newcommand{\Hmbar}{\overline{\Hm}}
\newcommand{\Delm}{\Delta_{\scriptscriptstyle\mbox{\rm\tiny M}}}
\newcommand{\Ax}{A_x}
\newcommand{\ax}{a_{\scriptscriptstyle x}}
\newcommand{\bx}{b_{\scriptscriptstyle x}}
\newcommand{\Az}{A_z}
\newcommand{\az}{a_{\scriptscriptstyle z}}
\newcommand{\bz}{b_{\scriptscriptstyle z}}
\newcommand{\Gbar}{\overline{\Gamma}}
\newcommand{\Wx}{W_x}
\newcommand{\Wz}{W_z}
\newcommand{\Wtilx}{\Wtil_x}
\newcommand{\Wtilz}{\Wtil_z}
\newcommand{\Qx}{Q_x}
\newcommand{\Qz}{Q_z}
\newcommand{\Qhatx}{\Qhat_x}
\newcommand{\Qhatz}{\Qhat_z}

\begin{document}

\title{A remark on manageable multiplicative unitaries\thanks{Partially
supported by Komitet
Bada\'{n} Naukowych, grant No. 2 P0A3 030 14 and the Foundation for Polish
Science}}
\author{P. M. So\l{}tan\qquad\qquad S. L. Woronowicz\\
Department of Mathematical Methods in Physics\\
Faculty of Physics, University of Warsaw\\
Ho\.{z}a 74, 00-682 Warsaw, Poland
}
\date{}
\maketitle

\vspace{-7mm}\begin{abstract}
We propose a weaker condition for multiplicative unitary operators related to
quantum groups, than the condition of manageability introduced by S.L.
Woronowicz.
We prove that all the main results of the theory of manageable multiplicative
unitaries
remain true under this weaker condition. We also show that multiplicative
unitaries arising
naturally in the construction of some recent examples of non-compact quantum
groups
satisfy our condition, but fail to be manageable.
\end{abstract}

\Section{Introduction}\label{intro}
The theory of multiplicative unitary operators initiated by S. Baaj and G.
Skandalis in \cite{baaj} has played a central role in the modern approach to
quantum groups. A unitary operator $W\in B(H\tens H)$ is called
multiplicative if it satisfies the pentagon equation (cf. \cite{baaj}):
\begin{equation}\label{Penta}
W_{23}W_{12}=W_{12}W_{13}W_{23}.
\end{equation}
However this condition alone does not guarantee that $W$ is a multiplicative
unitary related to a quantum group. S. Baaj and G. Skandalis proposed a
condition called regularity which unfortunately did not fit all applications
(cf. \cite{baajE2}, Proposition 4.2). In \cite{mu} the condition of regularity
was replaced by another one  called manageability. In \cite{mnw} it is shown
that all quantum groups possess a manageable multiplicative unitary which
is called the Kac-Takesaki operator.

As one might expect the manageability condition is often difficult to check
in particular examples. Moreover the natural choice for the multiplicative
unitaries in specific examples like the quantum ``$ax+b$'' and ``$az+b$''
groups turns out not to be manageable (cf. \cite{azb}, \cite{axb} and Section
\ref{appl}).

The aim of this paper is to weaken the manageability condition in such a way
that it suits the above mentioned examples (cf. Section \ref{appl}). The
condition we propose is the following: let $H$ be a Hilbert space and let
$W\in B(H\tens H)$ be a multiplicative unitary. We will suppose that there
exist two positive selfadjoint operators $\Qhat$ and $Q$ on $H$ with
$\ker{Q}=\ker{\Qhat}=\{0\}$ and a  unitary $\Wtil\in B(\Hbar\tens H)$ such
that
\[
W^*\bigl(\Qhat\tens Q\bigr)W=\Qhat\tens Q
\]
and
\[
\its{x\tens u}{W}{z\tens y}=
\its{\zbar\tens Qu}{\Wtil}{\xbar\tens Q^{-1}y}
\]
for all $x,z\in H$, $y\in D(Q^{-1})$ and $u\in D(Q)$. We hereby take
opportunity to change the name ``manageable'' and we shall call a
multiplicative unitary satisfying the above condition a modular multiplicative
unitary. The modularity is reflected in the existence of the scaling group and
the polar decomposition of the coinverse (cf. Theorem \ref{main}).

In case $\Qhat=Q$ we retain manageability and in particular any manageable
multiplicative unitary is modular. It turns out that all the results obtained
in \cite{mu} are true with these weaker assumptions. One may try to adapt the
proofs from \cite{mu} to this new situation, however we encountered some
difficulties with this programme. Instead we will construct a new
multiplicative unitary (on a different Hilbert space) which is manageable and
describes the same quantum group.

In Section \ref{mod} we will use an auxiliary separable Hilbert space $K$ and
a pair $(r,s)$ of closed operators acting on $K$ such that $s$ is selfadjoint,
$r$ is positive selfadjoint and
\[
r^{it}sr^{-it}=s-tI
\]
for all $t\in\RR$. An example of such a pair $(r,s)$ on $K=L^2(\RR)$ can be
obtained by taking
\[
\bigl(sf\bigr)(x)=xf(x),\qquad f\in L^2(\RR),\quad x\in\RR
\]
and letting $r$ be the analytic generator of the translation group:
\[
\bigl(r^{it}f\bigr)(x)=f(x-t),\qquad f\in L^2(\RR),\quad x,t\in\RR
\]
i.e. $r=\exp(-D)$ where $D=\frac{1}{i}\partial_x$.

Let us briefly recall the leg numbering notation which we already used in
\refeq{Penta}. Suppose $H$ is a Hilbert space and $T$ is an operator in $H$.
Then by $T_k$ we shall denote the operator
\[
\underbrace{I\tens\cdots\tens I}_{k-1}
\tens T\tens \underbrace{I\tens\cdots\tens I}_{n-k}
\]
acting on $H^{\tens n}$. A more sophisticated version of this notational
convention applies to operators acting on a tensor product of $H$ with itself.
Let $U\in B(H\tens H)$. Then $U_{kl}$ denotes the operator acting as $U$
on the $k$-th and $l$-th copies of $H$ sitting inside $H^{\tens n}$ and as
identity on all remaining copies of $H$ in $H^{\tens n}$. We say that this
operator has legs in the $k$-th and $l$-th factors of the tensor product
$H^{\tens n}$. We will also be using this notation when dealing with tensor
products of different Hilbert spaces.

Let $H$ be a separable Hilbert space and let $\Hbar$ be the complex conjugate
of $H$. For any $x\in H$ the corresponding element of $\Hbar$ will be denoted
by $\xbar$. Then $H\ni x\mapsto\xbar\in\Hbar$ is an antiunitary map. In
particular $\is{\xbar}{\overline{y}}=\is{y}{x}$ for any $x,y\in H$. For any
closed operator $m$ acting on $H$ the symbol $m^\top$ will denote the
transpose of $m$. By definition $D(m^\top)=\overline{D(m^*)}$ and
\[
m^\top\xbar=\overline{m^*x}
\]
for any $x\in D(m^*)$. If $m\in B(H)$ then $m^\top$ is the bounded operator on
$\Hbar$ such that $\its{\xbar}{m^\top}{\overline{y}}=\its{y}{m}{x}$ for all
$x,y\in H$. Clearly $B(H)\ni m\mapsto m^\top\in B(\Hbar)$ is an
antiisomorphism of $C^*$-algebras. Setting $\overline{\xbar}=x$ we identify
$\overline{\Hbar}$ with $H$. With this identification $m^{\top\top}=m$ for any
$m\in B(H)$.

\Section{The results}

\begin{Def}\label{bmmu}
Let $H$ be a Hilbert space and let $W\in B(H\tens H)$ be a multiplicative
unitary operator. We say that $W$ is modular if there exist two positive
selfadjoint operators $Q$ and $\Qhat$ on $H$ and a unitary operator
$\Wtil\in B(\Hbar\tens H)$ such that $\ker{Q}=\ker{\Qhat}=\{0\}$,
\begin{equation}\label{bmmu1}
W\bigl(\Qhat\tens Q\bigr)W^*=\Qhat\tens Q
\end{equation}
and
\begin{equation}\label{bmmu2}
\its{x\tens u}{W}{z\tens y}=\its{\zbar\tens Qu}{\Wtil}{\xbar\tens Q^{-1}y}
\end{equation}
for all $x,z\in H$, $u\in D(Q)$ and $y\in D(Q^{-1})$.
\end{Def}

We begin with an analogue of Proposition 1.4 of \cite{mu}. It shows that
the dual multiplicative unitary (cf. \cite{baaj}) of a modular multiplicative
unitary is modular. The operators $Q$ and $\Qhat$ exchange their positions.

\begin{Prop}
Let $H$ be a Hilbert space, $W$ a modular multiplicative unitary and $Q,\Qhat$
and $\Wtil$ the operators related to $W$ in the way described in Definition
\ref{bmmu}. Then
\begin{enumerate}
\item $\Wtil$and $\Qhat^{\top}\tens Q^{-1}$ commute.

\item
For any $x\in D(\Qhat^{-1})$, $z\in D(\Qhat)$ and $y,u\in H$ we have
\begin{equation}\label{contin}
\its{x\tens u}{W}{z\tens y}=\its{\overline{\Qhat z}\tens u}{\Wtil}
{\overline{\Qhat^{-1}x}\tens y}.
\end{equation}

\item The multiplicative unitary $\What=\Sigma W^*\Sigma$ is modular.
\end{enumerate}
\end{Prop}

\proof
The proof is almost the same as that of Proposition 1.4 of \cite{mu}. The
necessary modifications are so easy that we present only the proof of
Statement 3 as an example. It is obvious that $\What$ commutes with
$Q\tens\Qhat$. Moreover introducing the unitary
$\Whattil=\left(\Sigma\Wtil^*\Sigma\right)^{\top\tens\top}$ we have:
\begin{equation}\label{dualmod}
\its{x\tens u}{\What}{z\tens y}=
\its{\zbar\tens\Qhat u}{\Whattil}{\xbar\tens\Qhat^{-1}y}
\end{equation}
for any $x,z\in H$, $u\in D(\Qhat)$ and $y\in D(\Qhat^{-1})$. Indeed: using
in the fourth step formula \refeq{contin} we obtain
\[
\begin{array}{c}
\its{\zbar\tens\Qhat u}{\left(\Sigma\Wtil^*\Sigma\right)^{\top\tens\top}}
{\xbar\tens\Qhat^{-1}y}
=
\its{x\tens\overline{\Qhat^{-1}y}}
{\Sigma\Wtil^*\Sigma}{z\tens\overline{\Qhat u}}\vs\\
=
\its{\overline{\Qhat^{-1}y}\tens x}{\Wtil^*}{\overline{\Qhat u}\tens z}
=\overline{
\its{\overline{\Qhat u}\tens z}{\Wtil}{\overline{\Qhat^{-1}y}\tens x}}
=
\overline{\its{y\tens z}{W}{u\tens x}}\vs\\
=\its{u\tens x}{W^*}{y\tens z}=\its{x\tens u}{\Sigma W^*\Sigma}{z\tens y}
\end{array}
\]
and \refeq{dualmod} follows. It shows that $\What$ is modular.
\qed

Now we will present the main result of the paper.

\begin{Thm}\label{main}
Let $H$ be a separable Hilbert space and let $W\in B(H\tens H)$ be a modular
multiplicative unitary. Define
\begin{equation}\label{AAhat}
\left.\begin{array}{rcl}
A&=&{\left\{(\omega\tens\id)W:\omega\in B(H)_*\Vs{3.6}\right\}}^{\nc},\vs\\
\Ahat&=&{\left\{(\id\tens \omega)W^*:\omega\in B(H)_*\Vs{3.6}\right\}}^{\nc}.
\end{array}\right\}
\end{equation}
Then
\begin{enumerate}
\item\label{Calgs} $A$ and $\Ahat$ are nondegenerate separable
$C^*$-subalgebras in $B(H)$.\vsn

\item\label{WinM} $W\in M\bigl(\Ahat\tens A\bigr)$.\vsn

\item\label{comult} There exists a unique $\Del\in\Mor(A,A\tens A)$ such that
\begin{equation}\label{DelW}
(\id\tens\Del)W=W_{12}W_{13}.
\end{equation}
Moreover

\vsnp\hspace{1cm}\parbox{13.5cm}{\begin{trivlist}\sloppy
\item[{\rm(}i\/{\rm)}] $\Del$ is coassociative:
$(\Del\tens\id)\circ\Del=(\id\tens\Del)\circ\Del$,\vsnn
\item[{\rm(}ii\/{\rm)}]\sloppy
$\left\{\Del(a)(I\tens b):a,b\in A\Vs{3.6}\right\}$
and $\left\{(a\tens I)\Del(b):a,b\in A\Vs{3.6}\right\}$ are linearly
dense subsets of $A\tens A$.
\end{trivlist}}\vsnpp

\item\label{coinv}\sloppy There exists a unique closed linear operator $\ka$ on
the
Banach space $A$ such that \hbox{$\{(\omega\tens\id)W:\omega\in B(H)_*\}$} is a
core
for $\ka$ and
\[
\ka\bigl((\omega\tens\id)W\bigr)=(\omega\tens\id)W^*
\]
for any $\omega\in B(H)_*$. Moreover

\vsnp\hspace{1cm}\parbox{13.5cm}{\begin{trivlist}
\item[{\rm (}i\/{\rm)}] the domain of $\ka$ is a subalgebra of $A$ and $\ka$ is
antimultiplicative: for any $a,b\in D(\ka)$ we have $ab\in D(\ka)$ and
$\ka(ab)=\ka(b)\ka(a)$,\vsnn
\item[{\rm(}ii\/{\rm)}] \sloppy
the image $\ka(D(\ka))$ coincides with $D(\ka)^*$ and
$\ka(\ka(a)^*)^*=a$ for all \hbox{$a\in D(\ka)$},\vsnn
\item [{\rm(}iii\/{\rm)}] the operator $\ka$ admits the following polar
decomposition:
\[
\ka=R\circ\tau_{i/2},
\]
where $\tau_{i/2}$ is the analytic generator of a one parameter group
$\{\tau_t\}_{t\in\RR}$ of $*$-automorphisms of the $C^*$-algebra $A$ and $R$
is an involutive normal antiautomorphism of $A$,\vsnn
\item [{\rm(}iv\/{\rm)}] \sloppy
$R$ commutes with automorphisms $\tau_t$ for all
$t\in\RR$, in particular \hbox{$D(\ka)=D(\tau_{i/2})$},\vsnn
\item [{\rm(}v\/{\rm)}] $R$ and $\{\tau_t\}_{t\in\RR}$ are uniquely
determined.
\end{trivlist}}\vsnpp

\item\label{wzory} We have

\vsnp\hspace{1cm}\parbox{13.5cm}{\begin{trivlist}
\item[{\rm(}i\/{\rm)}]
$\Del\circ\tau_t=(\tau_t\tens\tau_t)\circ\Del$ for all $t\in\RR$,\vsnn
\item[{\rm(}ii\/{\rm)}]
$\Del\circ R=\sigma(R\tens R)\Del$,
\end{trivlist}}\vsnpp

where $\sigma$ denotes the flip map
$\sigma\colon A\tens A\ni a\tens b\mapsto b\tens a\in A\tens A$.\vsn

\item\label{QW} Let $\Wtil$ and $Q$ be the operators related to $W$ as in
Definition \ref{bmmu}. Then

\vsnp\hspace{1cm}\parbox{13.5cm}{\begin{trivlist}
\item[{\rm(}i\/{\rm)}] for any $t\in\RR$ and $a\in A$ we have
$\tau_t(a)=Q^{2it}aQ^{-2it}$,\vsnn
\item[{\rm(}ii\/{\rm)}] writing $a^R$ instead of $R(a)$ we have
$W^{\top\tens R}=\Wtil^*$.
\end{trivlist}}\vsnpp
\end{enumerate}
\end{Thm}

Apart from Statement \ref{wzory} the conclusion of Theorem \ref{main} is
the same as that of Theorem 1.5 of \cite{mu}. The only difference lies in the
weaker condition imposed on $W$.

\Section{The modified multiplicative unitary}\label{mod}
In this section for a given modular multiplicative unitary $W$ acting on
$H\tens H$ we shall construct unitary $\Wm$ acting on $\Hm\tens\Hm$. The
Hilbert space $\Hm=K\tens H$, where $K$ is the Hilbert space with a pair
$(r,s)$ of operators described in Section \ref{intro}.

Let $W\in B(H\tens H)$ be a modular multiplicative unitary. Define a unitary
operator \hbox{$X\in B(\Hm)$} by
\begin{equation}\label{X}
X=(I\tens Q)^{i(s\tens I)}(I\tens\Qhat)^{-i(s\tens I)}
=Q_2^{is_1}\Qhat_2^{-is_1}.
\end{equation}

Let us notice that
\begin{equation}\label{tozs}
X^*(r\tens Q)X=r\tens\Qhat.
\end{equation}
Indeed: for any $t\in\RR$ we have
\[
\begin{array}{rcl}
(r\tens Q)^{it}X
&=&(r^{it}\tens Q^{it})X=r_1^{it}Q_2^{it}Q_2^{is_1}\Qhat_2^{-is_1}\vs\\
&=&r_1^{it}Q_2^{i(s_1+tI)}\Qhat_2^{-is_1}
=Q_2^{is_1}r_1^{it}\Qhat_2^{-is_1}\vs\\
&=&Q_2^{is_1}\Qhat_2^{-i(s_1-tI)}r_1^{it}
=Q_2^{is_1}\Qhat_2^{-is_1}\Qhat_2^{it}r_1^{it}=X(r\tens\Qhat)^{it}
\end{array}
\]
and \refeq{tozs} follows.
Using the same method one can easily check that \refeq{bmmu1} implies that
\begin{equation}\label{trick}
Q_2^{it}WQ_2^{-it}=\Qhat_1^{-it}W\Qhat_1^{it}.
\end{equation}

Now we can define a unitary operator $\Wm\in B(K\tens H\tens K\tens H)$:
\begin{equation}\label{Wm}
\Wm=X_{12}W_{24}X_{12}^*.
\end{equation}
Notice that
\begin{equation}\label{albeW}
\Wm=(\alpha\tens\beta)W,
\end{equation}
where $\alpha$ and $\beta$ are injective unital and normal $*$-homomorphisms
\begin{equation}\label{defalbe}
\left.\begin{array}{rcl}
\alpha\colon B(H)\ni m&\longmapsto&X(I\tens m)X^*\in B(K\tens H),\vs\\
\beta\colon B(H)\ni m&\longmapsto&I\tens m\in B(K\tens H).
\end{array}\right\}
\end{equation}

\begin{Prop}
$\Wm$ is a manageable multiplicative unitary acting on $\Hm\tens\Hm$.
\end{Prop}

\proof
First we shall prove that $\Wm$ is a multiplicative unitary. We have to verify
that
\[
(\Wm)_{23}(\Wm)_{12}=(\Wm)_{12}(\Wm)_{13}(\Wm)_{23}
\]
which reads as
\begin{equation}\label{veri}
X_{34}W_{46}X_{34}^*X_{12}W_{24}X_{12}^*
=X_{12}W_{24}X_{12}^*X_{12}W_{26}X_{12}^*X_{34}W_{46}X_{34}^*
\end{equation}
on $K\tens H\tens K\tens H\tens K\tens H$. By commuting $X_{12}$ through
$X_{34}W_{46}X_{34}^*$ on the left hand side of \refeq{veri} and moving
$X_{12}^*$ though $X_{34}W_{46}X_{34}^*$ on the right hand side of
\refeq{veri} one reduces \refeq{veri} to
\begin{equation}\label{redu}
X_{34}W_{46}X_{34}^*W_{24}=W_{24}W_{26}X_{34}W_{46}X_{34}^*.
\end{equation}

The pentagon equation \refeq{Penta} gives us
\begin{equation}\label{dblpenta}
W_{46}W_{24}=W_{24}W_{26}W_{46}
\end{equation}
Taking into account \refeq{trick} and using the fact that operators with
different
legs commute we infer
that the right hand side of \refeq{dblpenta} is equal to
\[
\begin{array}{rcl}
{\rm RHS}&=&W_{24}\Qhat_2^{-is_3}
\bigl(\Qhat_2^{is_3}W_{26}\Qhat_2^{-is_3}\bigr)\Qhat_2^{is_3}W_{46}\vs\\
&=&W_{24}\Qhat_2^{-is_3}
\bigl(Q_6^{-is_3}W_{26}Q_6^{is_3}\bigr)\Qhat_2^{is_3}W_{46}\vs\\
&=&Q_6^{-is_3}W_{24}\Qhat_2^{-is_3}W_{26}Q_6^{is_3}W_{46}\Qhat_2^{is_3}.
\end{array}
\]
Thus
\begin{equation}\label{step2}
W_{46}W_{24}
=Q_6^{-is_3}W_{24}\Qhat_2^{-is_3}W_{26}Q_6^{is_3}W_{46}\Qhat_2^{is_3}.
\end{equation}
Applying the map
$m\mapsto Q_6^{is_3}\Qhat_2^{is_3}m\Qhat_2^{-is_3}Q_6^{-is_3}$ to both sides
of \refeq{step2} and repeatedly using \refeq{trick}  we obtain
\[
\begin{array}{rcl}
Q_6^{is_3}\Qhat_2^{is_3}W_{46}W_{24}\Qhat_2^{-is_3}Q_6^{-is_3}
&=&\Qhat_2^{is_3}W_{24}\Qhat_2^{-is_3}W_{26}
Q_6^{is_3}W_{46}Q_6^{-is_3}\vs\\
\bigl(Q_6^{is_3}W_{46}Q_6^{-is_3}\bigr)
\bigl(\Qhat_2^{is_3}W_{24}\Qhat_2^{-is_3}\bigr)
&=&\bigl(\Qhat_2^{is_3}W_{24}\Qhat_2^{-is_3}\bigr)W_{26}
\bigl(Q_6^{is_3}W_{46}Q_6^{-is_3}\bigr)\vs\\
\bigl(\Qhat_4^{-is_3}W_{46}\Qhat_4^{is_3}\bigr)
\bigl(Q_4^{-is_3}W_{24}Q_4^{is_3}\bigr)
&=&\bigl(Q_4^{-is_3}W_{24}Q_4^{is_3}\bigr)W_{26}
\bigl(\Qhat_4^{-is_3}W_{46}\Qhat_4^{is_3}\bigr).
\end{array}
\]
Now we apply the map $m\mapsto Q_4^{is_3}mQ_4^{-is_3}$ to both sides of the
last equality and use \refeq{X} to obtain \refeq{redu} which proves that
$\Wm$ is a
multiplicative unitary.

In order to prove manageability of $\Wm$ we have to construct the operators
required by Definition 1.2 of \cite{mu}. Let
\begin{equation}\label{Qm}
\Qm=r\tens Q
\end{equation}
and
\[
\Wmtil=\left(X_{12}^{\top}\right)^*\Wtil_{24}X_{12}^{\top},
\]
where the symbol $\top$ denotes the transposition
$B(\Hm)\ni m\longmapsto m^{\top}\in B(\Hmbar)$
(cf. Section \ref{intro}). Take $\xi,\xi'\in K\tens H$ and
$\eta,\eta'\in D(r)\tens_{\mbox{\rm\tiny alg}}D(Q)\subset D(\Qm)$.
Using selfadjointness of $r$, the equation
\refeq{bmmu} and the fact that operators with different legs commute we
obtain:
\[
\begin{array}{l}
\hspace{-7mm}\its{\xi\tens\eta}{\Wm}{\xi'\tens\Qm\eta'}
=\its{\xi\tens\eta}{X_{12}W_{24}X_{12}^*}{\xi'\tens\Qm\eta'}\vs\\
=\its{X^*\xi\tens\eta}{W_{24}}{X^*\xi'\tens\Qm\eta'}
=\its{X^*\xi\tens\eta}{W_{24}}{X^*\xi'\tens r_1Q_2\eta'}\vs\\
=\its{X^*\xi\tens\eta}{W_{24}}{r_3(X^*\xi'\tens Q_2\eta')}
=\its{r_3(X^*\xi\tens\eta)}{W_{24}}{X^*\xi'\tens Q_2\eta'}\vs\\
=\its{X^*\xi\tens r_1\eta}{W_{24}}{X^*\xi'\tens Q_2\eta'}
=\its{\overline{X^*\xi'}\tens Q_2r_1\eta}{\Wtil_{24}}
{\overline{X^*\xi}\tens\eta'}\vs\\
=\its{X^{\top}\overline{\xi'}\tens\Qm\eta}{\Wtil_{24}}
{X^{\top}\overline{\xi}\tens\eta'}
=\its{\overline{\xi'}\tens\Qm\eta}
{\left(X_{12}^\top\right)^*\Wtil_{24}X_{12}^\top}
{\overline{\xi}\tens\eta'}\vs\\
=\its{\overline{\xi'}\tens\Qm\eta}{\Wmtil}{\overline{\xi}\tens\eta'}.
\end{array}
\]
Since $D(r)\tens_{\mbox{\rm\tiny alg}}D(Q)$ is a core for $\Qm$ we have
\[
\its{\xi\tens\eta}{\Wm}{\xi'\tens\Qm\eta'}
=\its{\overline{\xi'}\tens\Qm\eta}{\Wmtil}{\overline{\xi}\tens\eta}
\]
for all $\eta,\eta'\in D(\Qm)$. Now replacing $\eta'$ by
$\Qm^{-1}\eta'$ we get
\[
\its{\xi\tens\eta}{\Wm}{\xi'\tens\eta'}
=\its{\overline{\xi'}\tens\Qm\eta}{\Wmtil}{\overline{\xi}\tens\Qm^{-1}\eta'}
\]
for all $\xi,\xi'\in K\tens H$, $\eta\in D(\Qm)$ and $\eta'\in D(\Qm^{-1})$.

\sloppy
It remains to prove that $\Wm$ commutes with $\Qm\tens\Qm$. Using formula
\refeq{tozs}, \refeq{bmmu1} and again \refeq{tozs} we obtain
\[
\begin{array}{rcccl}
(\Qm\tens\Qm)\Wm
&=&(r\tens Q\tens r\tens Q)(X\tens I\tens I)W_{24}(X^*\tens I\tens I)\vs\\
&=&(X\tens I\tens I)(r\tens\Qhat\tens r\tens Q)W_{24}(X^*\tens I\tens I)\vs\\
&=&(X\tens I\tens I)W_{24}(r\tens\Qhat\tens r\tens Q)(X^*\tens I\tens I)\vs\\
&=&(X\tens I\tens I)W_{24}(X^*\tens I\tens I)(r\tens Q\tens r\tens Q)
&=&\Wm(\Qm\tens\Qm).
\end{array}
\]
We have thus checked that $\Wm$ satisfies all the conditions of Definition
1.2 from \cite{mu}.
\qed

We are now free to use the theory presented in \cite{mu}. All objects
constructed for $\Wm$ with help of Theorem 1.5 of that paper will be denoted
by letters with a subscript {\tiny M}. For example
\begin{equation}\label{AmiAm}
\left.\begin{array}{rcl}
\Am&=&\left\{(\Phi\tens\id)\Wm:\Phi\in B(\Hm)_*\Vs{3.6}\right\}^{\nc},\vs\\
\Amhat&=&\left\{(\id\tens \Phi)\Wm^*:\Phi\in B(\Hm)_*\Vs{3.6}\right\}^{\nc}.
\end{array}\right\}
\end{equation}
We also have $\Delm$, $\kam$, $\{\taum_t\}_{t\in\RR}$ and $\Rm$.

\Section{Proof of Theorem \ref{main}}
{\bf Ad \ref{Calgs}.} \sloppy
We know that $\Am$ and $\Amhat$ defined by \refeq{AmiAm} are nondegenerate
separable \hbox{$C^*$-subalgebras} of $B(\Hm)$. Recall that $\alpha$
and $\beta$ are ultra-weakly continuous injections of $B(H)$ into $B(\Hm)$.
Therefore for any normal functional $\omega$ on $B(H)$ there exits
$\Phi,\Phi'\in B(\Hm)_*$ such that
\[
\omega=\Phi\circ\alpha=\Phi'\circ\beta.
\]
Keeping this fact in mind, remembering the definitions \refeq{AAhat} and
formula \refeq{albeW} we have
\[
\begin{array}{rcl}
\beta(A)&=&\left\{\beta\bigl((\omega\tens\id)W\bigr):\omega\in B(H)_*\Vs{3.6}
\right\}^{\nc}\vs\\
&=&\left\{(\omega\tens\beta)W:\omega\in B(H)_*\Vs{3.6}\right\}^{\nc}\vs\\
&=&\left\{(\Phi\circ\alpha\tens\beta)W:\Phi\in B(\Hm)_*\Vs{3.6}\right\}^{\nc}
\vs\\
&=&\left\{(\Phi\tens\id)\Wm:\Phi\in B(\Hm)_*\Vs{3.6}\right\}^{\nc}=\Am.
\end{array}
\]
Similarly we prove that
\[
\alpha\bigl(\Ahat\bigr)=\Amhat.
\]
Now it is easy to see that $A$ and $\Ahat$ are nondegenerate separable
$C^*$-subalgebras in $B(H)$.

{\bf Ad \ref{WinM}.} We know that $\Wm\in M\bigl(\Amhat\tens\Am\bigr)$. In
other words (cf. \refeq{albeW})
\[
(\alpha\tens\beta)W\in M\bigl(\alpha\bigl(\Ahat\bigr)\tens\beta(A)\bigr)=
M\bigl((\alpha\tens\beta)\bigl(\Ahat\tens A\bigr)\bigr)
\]
and it follows that $W\in M\bigl(\Ahat\tens A\bigr)$.

{\bf Ad \ref{comult}.} We have $\Delm\in\Mor(\Am,\Am\tens\Am)$ and a
$*$-isomorphism $\beta\colon A\to\Am$, so we can define
$\Del=(\beta\tens\beta)^{-1}\Delm\beta$. This provides us with a coassociative
$\Del\in\Mor(A,A\tens A)$ such that
$\left\{\Del(a)(I\tens b):a,b\in A\Vs{3.6}\right\}$ and
$\left\{(a\tens I)\Del(b):a,b\in A\Vs{3.6}\right\}$ are linearly dense subsets
of $A\tens A$. Furthermore notice that
\[
\begin{array}{rcl}
(\id\tens\Del)W&=&\bigl(\id\tens(\beta\tens\beta)^{-1}\Delm\beta\bigr)W\vs\\
&=&\bigl(\alpha^{-1}\alpha\tens(\beta\tens\beta)^{-1}\Delm\beta\bigr)W\vs\\
&=&\bigl(\alpha^{-1}\tens(\beta\tens\beta)^{-1}\Delm\bigr)
(\alpha\tens\beta)W\vs\\
&=&\bigl(\alpha^{-1}\tens(\beta\tens\beta)^{-1}\Delm\bigr)\Wm\vs\\
&=&(\alpha\tens\beta\tens\beta)^{-1}(\id\tens\Delm)\Wm\vs\\
&=&(\alpha\tens\beta\tens\beta)^{-1}(\Wm)_{12}(\Wm)_{13}\vs\\
&=&(\alpha\tens\beta\tens\beta)^{-1}(\Wm)_{12}
(\alpha\tens\beta\tens\beta)^{-1}(\Wm)_{13}\vs\\
&=&(\alpha\tens\beta\tens\id)^{-1}(\Wm)_{12}
(\alpha\tens\id\tens\beta)^{-1}(\Wm)_{13}=W_{12}W_{13}.
\end{array}
\]

{\bf Remark.}\quad
Despite a fairly complicated way of introducing the comultiplication on $A$ we
can still recover formula (5.1) of \cite{mu}
i.e.
\begin{equation}\label{DelRem}
\Del(a)=W(a\tens I)W^*
\end{equation}
for all $a\in A$ (cf. \cite{baaj}, Th\'{e}or\`{e}me 3.8). Indeed: take
$a=(\omega\tens\id)W$ then using \refeq{DelW} and \refeq{Penta} we obtain
\[
\begin{array}{rcl}
\Del(a)
&=&(\omega\tens\id\tens\id)(\id\tens\Del)W\vs\\
&=&(\omega\tens\id\tens\id)W_{12}W_{13}\vs\\
&=&(\omega\tens\id\tens\id)W_{12}W_{13}\vs\\
&=&(\omega\tens\id\tens\id)W_{23}W_{12}W_{23}^*\vs\\
&=&W\bigl((\omega\tens\id)W\tens I\bigr)W^*.
\end{array}
\]
For an arbitrary $a\in A$ we use the continuity argument. This also proves
the uniqueness of $\Del$.

{\bf Ad \ref{coinv}.} Since $\beta$ is an isomorphism $A\to\Am$ we can define
$\ka=\beta^{-1}\kam\beta$.  Now it is important to notice (cf. the proof of
Statement \ref{Calgs}) that
\begin{equation}\label{phial}
\left.\begin{array}{rcl}
\beta^{-1}\bigl((\Phi\tens\id)\Wm\bigr)&=&(\Phi\circ\alpha\tens\id)W\vs\\
\beta^{-1}\bigl((\Phi\tens\id)\Wm^*\bigr)&=&(\Phi\circ\alpha\tens\id)W^*.
\end{array}\right\}
\end{equation}
Then first of all it follows from \refeq{phial} that
\[
\beta\left(\left\{(\omega\tens\id)W:\omega\in B(H)_*\Vs{3.6}\right\}\right)
=\left\{(\Phi\tens\id)\Wm:\Phi\in B(\Hm)_*\Vs{3.6}\right\}.
\]
Furthermore for $\omega\in B(H)_*$ we have $\omega=\Phi\circ\alpha$ for some
$\Phi\in B(\Hm)_*$ and using \refeq{phial} we get
\[
\begin{array}{rcl}
\ka\bigl((\omega\tens\id)W\bigr)
&=&\beta^{-1}\kam\beta\bigl((\omega\tens\id)W\bigr)
=\beta^{-1}\kam\bigl((\Phi\tens\id)\Wm\bigr)\vs\\
&=&\beta^{-1}\bigl((\Phi\tens\id)\Wm^*\bigr)=(\omega\tens\id)W^*.
\end{array}
\]
Since $\left\{(\Phi\tens\id)\Wm:\Phi\in B(\Hm)_*\Vs{3.6}\right\}$ is a
core for $\kam$ we see that
$\left\{(\omega\tens\id)W:\omega\in B(H)_*\Vs{3.6}\right\}$ is a core for
$\ka$. Now setting
\begin{equation}\label{tauR}
\left.\begin{array}{rcl}
\tau_t&=&\beta^{-1}\taum_t\beta,\qquad t\in\RR,\vs\\
R&=&\beta^{-1}\Rm\beta
\end{array}\right\}
\end{equation}\sloppy
we see that assertions ({\em i}\/) -- ({\em v}\/) follow directly from
analogous statements for $\Am,\kam,\Rm$ and $\{\taum_t\}_{t\in\RR}$ (cf.
\cite{mu}, Theorem 1.5, Statement 4.) and the fact the $\beta$ is a normal
\hbox{$*$-isomorphism} of $A$ onto $\Am$.

{\bf Ad \ref{QW}.} We know (cf. \cite{mu}, Theorem 1.5, Statement 5) that
for any $\am\in\Am$ and any $t\in\RR$
\begin{equation}\label{impltau}
\taum_t(\am)=\Qm^{2it}\am\Qm^{-2it}.
\end{equation}
Thus formula ({\em i}\/) follows from the first line of \refeq{tauR},
\refeq{impltau}, \refeq{Qm} and the definition of $\beta$.

From the results of \cite{mu} (formula (1.14)) we know that
\begin{equation}\label{Wform}
\Wm^{\top\tens\Rm}=\Wmtil^*.
\end{equation}
Notice that
\begin{equation}\label{albeWt}
\Wmtil=(\alpha^\top\tens\beta)\Wtil,
\end{equation}
\sloppy where $\alpha^\top\colon B(\Hbar)\to B(\Hmbar)$ is a normal
\hbox{$*$-monomorphism} given by
\[
\alpha^\top(m)=\left(X^\top\right)^*(I\tens m)X^\top.
\]
It is easy to check that
\begin{equation}\label{altop}
\top\circ\alpha=\alpha^\top\circ\top.
\end{equation}
Finally recall that from the definition of $R$ \refeq{tauR} it follows that
\begin{equation}\label{betaR}
\Rm\circ\beta=\beta\circ R.
\end{equation}
Now taking into account \refeq{albeWt} and \refeq{Wform} and using
\refeq{albeW}, \refeq{altop} and \refeq{betaR} we obtain
\[
\begin{array}{rcl}
(\alpha^\top\tens\beta)\Wtil^*
&=&\Wmtil^*=\Wm^{\top\tens\Rm}\vs\\
&=&\bigl((\alpha\tens\beta)W\bigr)^{\top\tens\Rm}\vs\\
&=&(\alpha^\top\tens\beta)\left(W^{\top\tens R}\right),
\end{array}
\]
which gives formula ({\em ii}\/).

{\bf Ad \ref{wzory}.} Recall (cf. \cite{mu}, Theorem 1.5, Statement 5 and
formula (5.1)) that for any $\am\in\Am$ we have
\[
\Delm(\am)=\Wm(\am\tens I)\Wm^*.
\]
Now formula ({\em i}\/) follows from an easy computation:
\[
\begin{array}{rcl}
\Delm\bigl(\taum_t(\am)\bigr)
&=&\Wm\bigl(\taum_t(\am)\tens I\bigr)\Wm^*\vs\\
&=&\Wm(\Qm^{2it}\am\Qm^{-2it}\tens I)\Wm^*\vs\\
&=&\Wm(\Qm^{2it}\tens\Qm^{2it})(\am\tens I)
(\Qm^{-2it}\tens\Qm^{-2it})\Wm^*\vs\\
&=&(\Qm^{2it}\tens\Qm^{2it})\Wm(\am\tens I)
\Wm^*(\Qm^{-2it}\tens\Qm^{-2it})\vs\\
&=&(\taum_t\tens\taum_t)\Delm(\am),
\end{array}
\]
where in the second last equality we used the fact that $\Wm$ commutes with
$\Qm\tens\Qm$.


Consider a one parameter group
$\RR\ni t\mapsto\sigma_t\in{\rm Aut}\bigl(B(H)\bigr)$,
$\sigma_t(m)=\Qhat^{it}m\Qhat^{-it}$. Let $\sigma_{-i}$ be its analytic
generator (\cite{zsido}). Now for $m\in D(\sigma_{-i})$ define
\[
\Tau(m)=\left(\sigma_{-i}(m)\right)^\top.
\]
It follows (cf. \cite{zsido}) that $\Tau$ is a closed densely defined
injective operator whose domain is a subalgebra of $B(H)$ and that $\Tau$
is antimultiplicative.

\begin{Prop}\label{defTau}
For any $\omega\in B(H)_*$ the element $(\id\tens\omega)\in D(\Tau)$ and
\[
\Tau\bigl((\id\tens\omega)W\bigr)=(\id\tens\omega)\Wtil.
\]
\end{Prop}

\proof
From formula \refeq{contin} we infer that for any $\omega\in B(H)_*$
we have
\[
\its{x}{(\id\tens\omega)W}{z}
=\its{\overline{\Qhat z}}{(\id\tens\omega)\Wtil}{\overline{\Qhat^{-1}x}}
\]
for any $x\in D(\Qhat^{-1})$ and $z\in D(\Qhat)$. This can be rephrased as
\[
\its{\zbar}{(\id\tens\omega)\Wtil}{\xbar}=
\its{\Qhat x}{(\id\tens\omega)W}{\Qhat^{-1}z}
\]
for any $z\in D(\Qhat^{-1})$ and $x\in D(\Qhat)$. Therefore
$\Qhat\bigl((\id\tens\omega)W\bigr)\Qhat^{-1}$ extends to a bounded operator on
$H$
and
\[
\left(\Qhat\bigl((\id\tens\omega)W\bigr)\Qhat^{-1}\right)^\top
=(\id\tens\omega)\Wtil\in B(\Hbar).
\]
This shows that any element of the form $(\id\tens\omega)W$ lies in the domain
of $\Tau$ and that
\hbox{$\Tau\bigl((\id\tens\omega)W\bigr)=(\id\tens\omega)\Wtil$}.
\qed

Formula \refeq{DelRem} allows us to define $\Del(m)$ for any $m\in B(H)$
which justifies its use in the statement of the next proposition.

\begin{Prop} We have
\begin{equation}\label{DelWtil}
(\id\tens\Del)\Wtil=\Wtil_{13}\Wtil_{12}.
\end{equation}
\end{Prop}

\proof
Take $\mu,\nu\in B(H)_*$ and denote by $\nu*\mu$ the normal functional
$(\mu\tens\nu)\circ\Del$. Now using \refeq{DelW} and the fact that $\Tau$
defined in Proposition \ref{defTau} is antimultiplicative we compute:
\[
\begin{array}{rcl}
(\id\tens\mu\tens\nu)(\id\tens\Del)\Wtil
&=&(\id\tens\nu*\mu)\Wtil\vs\\
&=&\Tau\bigl((\id\tens\nu*\mu)W\bigr)\vs\\
&=&\Tau\bigl((\id\tens\mu\tens\nu)(\id\tens\Del)W\bigr)\vs\\
&=&\Tau\bigl((\id\tens\mu\tens\nu)W_{12}W_{13}\bigr)\vs\\
&=&\Tau\bigl((\id\tens\mu)W(\id\tens\nu)W\bigr)\vs\\
&=&\Tau\bigl((\id\tens\nu)W\bigr)\Tau\bigl((\id\tens\mu)W\bigr)\vs\\
&=&(\id\tens\nu)\Wtil(\id\tens\mu)\Wtil\vs\\
&=&(\id\tens\mu\tens\nu)\Wtil_{13}\Wtil_{12}.
\end{array}
\]
Now since functionals of the form $\mu\tens\nu$ separate elements of $B(H\tens
H)$
we obtain \refeq{DelWtil}.
\qed

We will use formula \refeq{DelWtil} to prove assertion ({\em ii}\/) of point
\ref{wzory} of our theorem. Applying $*$ to both sides of \refeq{DelWtil}
we get
\begin{equation}\label{DelWtilst}
(\id\tens\Del)\Wtil^*=\Wtil_{12}^*\Wtil_{13}^*.
\end{equation}
Notice that due to ({\em ii}\/) of point \ref{QW} of our theorem the left hand
side of \refeq{DelWtilst} is equal to
\[
(\top\tens\Del\circ R)W
\]
while the right hand side of \refeq{DelWtilst} equals
\[
\bigl(W^{\top\tens R}\bigr)_{12}\bigl(W^{\top\tens R}\bigr)_{13}
=(\top\tens R\tens R)W_{13}W_{12}
=(\top\tens R\tens R)(\id\tens\sigma)(\id\tens\Del)W.
\]
In other words
\begin{equation}\label{fin}
(\top\tens\Del\circ R)W
=(\top\tens R\tens R)(\id\tens\sigma)(\id\tens\Del)W
=(\top\tens\sigma\circ(R\tens R)\circ\Del)W.
\end{equation}
Applying $(\omega\circ\top\tens\id)$ to both sides of \refeq{fin} and taking
into
account \refeq{AAhat} we obtain formula ({\em ii}\/).

\Section{Applications}\label{appl}
In this section we will briefly present two examples of quantum groups whose
naturally
occurring multiplicative unitaries are modular, but not manageable. These groups
are the quantum
``$ax+b$'' and ``$az+b$'' groups constructed in \cite{axb} and \cite{azb}
respectively.
The algebras $\Ax$ and $\Az$ of continuous functions vanishing at infinity on
these groups
are generated (\cite{gen}) by $\ax,\ax^{-1},\bx,\beta$ and $\az,\az^{-1},\bz$
affiliated
with $\Ax$ and $\Az$. These distinguished elements are subject to relations
\[
\left(\begin{array}{c}
\ax\mbox{ and }\bx\mbox{ are selfadjoint}\\
\ax\mbox{ is strictly positive and}\\
\ax^{it}\bx\ax^{-it}=e^{\hbar t}\bx\\
\mbox{for any }t\in\RR,\\
\beta^2=\chi(\bx\neq0),\ \beta\ax=\ax\beta\\
\mbox{and }\beta\bx=-\bx\beta
\end{array}\right)
\qquad
\left(\begin{array}{c}
\az\mbox{ and }\bz\mbox{ are normal operators}\\
{\rm Sp\,}\az,{\rm Sp\,}\bz\subset\Gbar ,\ \ker{\az}=\{0\}\\
({\rm Phase\,}\az)\bz({\rm Phase\,}\az)^*=e^{\frac{2\pi i}{N}}\bz\\
|\az|^{-it}\bz|\az|^{it}=e^{\frac{2\pi}{N}t}\bz\\
\mbox{for any }t\in\RR
\end{array}\right)
\]
where $\Gbar=\{0\}\cup\bigcup\limits_{k=0}^{N-1}e^{\frac{2\pi i}{N}k}\RR_+$,
$\RR_+=\{x\in\RR:x>0\}$ and
$\hbar\in]0,\pi[$ and $N\in\NN$ are deformation parameters satisfying additional
requirements.
The natural choices for multiplicative unitary operators for these groups are
\begin{equation}\label{Wx}
\Wx=
F_{\hbar}\left(
e^{i\frac{i}{\hbar}}\bx^{-1}\ax\tens\bx,
ie^{\frac{i\pi^2}{2\hbar}}(\beta\tens\beta)\chi\left(\bx\tens\bx<0\right)\right)
^*
e^{
\frac{i}{\hbar}
\log{(|\bx|^{-1})}
\tens\log{\ax}
}
\end{equation}
for the ``$ax+b$'' group and
\begin{equation}\label{Wz}
\Wz=F_N\left(\az\bz^{-1}\tens\bz\right)\chi\left(\bz^{-1}\tens
I,I\tens\az\right)
\end{equation}
for the ``$az+b$'' group. In the above formulae we choose representations in
which $\bx$ and $\bz$ are invertible and $F_{\hbar}, F_N$ and $\chi$ are
special functions. It can be shown that  both $\Wx$ and $\Wz$ are
multiplicative unitary operators, but neither of them is manageable.
Nevertheless they are both modular with  \[
\begin{array}{rclcrcl}
\Qhatx&=&|\bx|^{\frac{1}{2}},&\quad&\Qx&=&(\ax)^{\frac{1}{2}},\vs\\
\Qhatz&=&|\bz|,&\quad&\Qz&=&|\az|
\end{array}
\]
and\sloppy
\[
\begin{array}{l}
\Wtilx=F_{\hbar}\left(
-e^{i\frac{\hbar}{2}}\left(\bx^{-1}\ax\right)^\top\tens
e^{i\frac{\hbar}{2}}\bx\ax^{-1},
-(\beta\tens\beta)\chi\bigl(e^{i\frac{\hbar}{2}}(\bx^{-1}\ax)^\top\tens\bx>0
\bigr)
\right)e^{\frac{i}{\hbar}\log{(\ax)^\top}\tens\log{\ax} },\vs\\
\Wtilz=F_N\left(-\left(\az\bz^{-1}\right)^\top\tens
e^{\frac{2\pi i}{N}}\az^{-1}\bz\right)^*\chi\left({\bz^{-1}}^\top\tens
I,I\tens\az\right ).
\end{array}
\]
In both constructions \cite{axb} and \cite{azb} a clever trick was used to
obtain manageability of the unitaries \refeq{Wx} and \refeq{Wz}. This trick
was the basis of our construction of the modified multiplicative unitary
presented in Section \ref{mod}. Using Theorem \ref{main} one is able to carry
out the construction of the two quantum groups without having to resort to
some slightly unintuitive means (cf. \cite{axb} Theorem 2.1 and \cite{azb}
Theorem 3.1).

\end{document}